\documentclass[a4paper, 11pt]{article}
\usepackage{amsthm}
\usepackage{amsfonts}
\usepackage{latexsym}
\usepackage{epic}
\usepackage{graphics}
\usepackage{graphicx}%
\usepackage{epsf, epstopdf}
\usepackage{etex}
\usepackage{pinlabel}

\usepackage[top=2.5cm,bottom=2.5cm,left=2.5cm,right=2.5cm]{geometry}
\usepackage{graphicx}
\usepackage{amsmath,amsthm,amssymb,lineno}

\usepackage{graphicx}           
\usepackage{color, ifpdf}  
\usepackage[dvipsnames]{xcolor}            
\usepackage{graphicx}
\usepackage{geometry}
\usepackage{float}
\usepackage{indentfirst}        
\usepackage{amsmath,amssymb,bm} 
\usepackage{cases}              
\usepackage{setspace}
\usepackage{graphicx,booktabs,multirow}
\usepackage{enumerate,mathdots}

\usepackage{latexsym}
\usepackage{graphicx,color, booktabs,multirow}
\usepackage{tikz,chemfig}
\usepackage{graphicx,booktabs,multirow,mathrsfs}
\usepackage{appendix}
\usepackage{yhmath}
\usepackage{setspace}

\newtheorem{theo}{Theorem}
\newtheorem{lemm}[theo]{Lemma}

\newtheorem{coro}[theo]{Corollary}

\newtheorem{clm}{Claim}

\newcommand \Clm[2]
{
\begin{clm}\label{#1}
#2
\end{clm}
}

\usetikzlibrary{backgrounds}
\usetikzlibrary{arrows}
\usetikzlibrary{shapes,shapes.geometric,shapes.misc}
\pgfdeclarelayer{edgelayer}
\pgfdeclarelayer{nodelayer}
\pgfsetlayers{background,edgelayer,nodelayer,main}
\tikzstyle{none}=[inner sep=0mm]
\tikzstyle{every loop}=[]

\tikzstyle{dotted}=[dash pattern=on \pgflinewidth off 2pt]
\tikzstyle{dashed}=[dash pattern=on 3pt off 3pt]

\newcommand \tikzp[2]
{
	\begin{center}
		\begin{tikzpicture}[scale=#1]
			#2
		\end{tikzpicture}
	\end{center}
}

\tikzstyle{new style 0}=[fill=black, draw=black, shape=circle]
\tikzstyle{red style 1}=[fill=red, draw=black, shape=circle]
\tikzstyle{blue style 2}=[fill=blue, draw=black, shape=circle]
\tikzstyle{white style 4}=[fill=white, draw=black, shape=circle]
\tikzstyle{bklack style 5}=[fill=black, draw=black, shape=rectangle]
\tikzstyle{red style 3}=[fill=red, draw=black, shape=rectangle]
\tikzstyle{yellow style 7}=[fill=yellow, draw=black, shape=rectangle]
\tikzstyle{new style 8}=[fill={rgb,255: red,0; green,132; blue,0}, draw={rgb,255: red,0; green,131; blue,0}, shape=circle]

\tikzstyle{new edge style 0}=[-]
\tikzstyle{new edge style 1}=[-, draw=red]
\tikzstyle{new edge style 2}=[-, draw=blue]
\tikzstyle{new edge style 3}=[-, draw={rgb,255: red,0; green,156; blue,0}]

\allowdisplaybreaks

\numberwithin{equation}{section}

\newcounter{countcase}

\newcounter{countclaim}

\def \proof {\noindent {\it Proof}. }
\newcommand{\claimend}{{\hfill$\natural $}}

\newcommand{\proofend}{
{\hfill$\Box$}
\setcounter{countclaim} {0}
\setcounter{countcase} {0}
}

\def \cov {{\mathcal {H}}}

\def \N {{\mathbb {N}}}

\newcommand \equ[2]
{
\begin{equation}\label{#1}
#2
\end{equation}
}

\newcommand \eqn[2]
{
\begin{eqnarray}\label{#1}
#2
\end{eqnarray}
}

\newcommand \aln[2]
{
\begin{align}\label{#1}
#2
\end{align}
}

\def \var {\frac{2k+2}{2k+1}}
\def \varn {2k+2}
\def \vard {2k+1}

\begin{document}
	
	\renewcommand{\contentsname}{\LARGE \center Content}         

\title{$Z_{DP}(n)$ is upperly bounded by $n^2-(n+3)/2$
}

\author{Meiqiao Zhang\thanks{Corresponding author. Email: nie21.zm@e.ntu.edu.sg and 
		meiqiaozhang95@163.com.} and Fengming Dong\thanks{Email: fengming.dong@nie.edu.sg (expired on 24/03/2027) and donggraph@163.com.}
	\\
	\small National Institute of Education,
	Nanyang Technological University, Singapore
}

\date{}
\maketitle

\begin{abstract} 
DP-coloring was introduced by Dvo\v{r}\'{a}k and Postle and is a generalization of proper coloring.
For any graph $G$,  
let $\chi(G)$ and $\chi_{DP}(G)$
denote the chromatic number and the DP-chromatic number of $G$
respectively. 
In this article, we show that 
$\chi_{DP}(G\vee K_s)=\chi(G\vee K_s)$ holds for 
$s=\left \lceil \frac{4(k+1)m}{2k+1} \right \rceil\le \lceil 2.4m\rceil$,
where $k=\chi(G)$, $m=|E(G)|$ and $G\vee K_s$ is the join of $G$ and 
the complete graph $K_s$.
Hence $Z_{DP}(n)\le n^2-(n+3)/2$ holds for every integer $n\ge 2$,
where $Z_{DP}(n)$ is the minimum natural number $s$ such that 
$\chi_{DP}(G\vee K_s)=\chi(G\vee K_s)$ holds for 
every graph $G$ of order $n$.
Our result improves the best current upper bound
$Z_{DP}(n)\le 1.5n^2$ due to 
Bernshteyn, Kostochka and Zhu.
\end{abstract}


\section{Introduction
\label{sec1}}

In this article, we consider simple graphs only. 
For any graph $G$, let $V(G)$ and $E(G)$ be its vertex set and edge set, respectively. 
For any vertex $u\in V(G)$, 
let $N_G(u)$ be the set of  vertices adjacent to $u$ in $G$,
let $d_G(u)$ (or simply $d(u)$) be the \textit{degree} of $u$ in $G$
(i.e., $d(u)=|N_G(u)|$), 
let $G_u$ be the subgraph obtained from $G$ by deleting $u$, 
and call $u$  {\it simplicial}
if either $d_G(u)=0$ or $N_G(u)$ is a clique.
For any non-empty subset $V_0$ of $V(G)$, 
let $G[V_0]$ denote the subgraph of $G$ induced by $V_0$, 
let $N_G(V_0)=\bigcup_{v\in V_0}N_G(v)$, and write $N_G(\{v_1,\cdots,v_\ell\})$ as $N_G(v_1,\cdots,v_\ell)$.
For any two disjoint subsets $V_1, V_2$ of $V(G)$,
let $E_G(V_1,V_2)$ be the set of edges $uv\in E(G)$, 
where $u\in V_1$ and $v\in V_2$. 
For any vertex-disjoint graphs $G_1$ and $G_2$, the {\it join} of $G_1$ and $G_2$,
denoted by $G_1\vee G_2$,  is the graph with vertex set $V(G_1)\cup V(G_2)$ 
and edge set $E(G_1)\cup E(G_2)\cup \{v_1v_2: v_1\in V(G_1)\text{ and }v_2\in V(G_2)\}$.

\subsection{Preliminary \label{sec1.1}}

Denote the set of non-negative integers by $\N$. For any positive integer $k$, 
let $[k]=\{1,\cdots,k\}$. 

For any graph $G$ and any positive integer $k$, a \textit{proper $k$-coloring} of $G$ is a mapping $f:V(G)\rightarrow[k]$, such that $f(u)\neq f(v)$ for all $uv\in E(G)$. A graph $G$ is \textit{$k$-colorable} if $G$ has a proper $k$-coloring. 
The \textit{chromatic number} $\chi(G)$ is the minimum positive integer $k$ 
such that $G$ is $k$-colorable. 
The chromatic number is one of the oldest and most important 
invariants in graph theory.

As a generalization of proper coloring, the notion of \textit{list coloring} was introduced by Vizing~\cite{Vizing1976} and Erd\H{o}s, Rubin and Taylor~\cite{Erdos1979} independently. For any graph $G$, a \textit{list assignment} of $G$ is a mapping $L:V(G)\rightarrow \mathcal{P}(\N)$, where $\mathcal{P}(\N)$ is the powerset of $\N$.  $L$ is a \textit{$k$-list assignment} of $G$ if $|L(v)|=k$ holds for all $v\in V(G)$. For any list assignment $L$ of $G$, a \textit{list coloring} is a mapping $f:V(G)\rightarrow \N$ with $f(v)\in L(v)$ for all $v\in V(G)$ and $f(u)\neq f(v)$ for all $uv\in E(G)$. The \textit{list chromatic number} $\chi_\ell(G)$ is the minimum positive integer $k$, such that $G$ has a list coloring for each $k$-list assignment $L$. 


To deal with a problem in list coloring, Dvo\v{r}\'{a}k and Postle~\cite{Dvorak2018} recently introduced the \textit{correspondence coloring}, which is referred as \textit{DP-coloring} more often afterwards. The formal definition is as follows.

For any graph $G$, a \textit{cover }of $G$ is an ordered pair $\cov=(L,H)$, where $H$ is a graph and $L$ is mapping from $V(G)$ to $\mathcal{P}(V(H))$ satisfying the conditions below:
\begin{enumerate}
\item the set $\{L(u):u\in V(G)\}$ is a partition of $V(H)$,
\item for every $u\in V(G)$, $H[L(u)]$ is a complete graph,
\item if $u$ and $v$ are not adjacent in $G$, then $E_H(L(u), L(v))=\emptyset$,
\item for each edge $uv\in E(G)$, $E_H(L(u), L(v))$ is a matching.
\end{enumerate}

A cover $\cov=(L,H)$ of $G$ is \textit{$m$-fold} 
if $|L(v)|=m$ for all $v\in V(G)$. 
For any cover $\cov=(L,H)$ of $G$,
an \textit{$\cov$-coloring} of $G$ is an independent set 
in $H$ of size $|V(G)|$ and $G$ is \textit{$\cov$-colorable} if $G$ has an $\cov$-coloring.
Obviously, any $\cov$-coloring $I$ of $G$ has the property that 
$|I\cap L(v)|=1$ for each $v\in V(G)$. 
The \textit{DP-chromatic number} $\chi_{DP}(G)$ is the minimum positive integer $m$ such that 
$G$ is $\cov$-colorable for each $m$-fold cover $\cov$ of $G$.

The {\it coloring number} of a graph $G$, denoted by $col(G)$, 
	is the smallest integer $r$ for which 
	there exists an ordering $v_1, v_2,\cdots, v_n$ of all the vertices in $G$ 
	such that $|N_G(v_i)\cap \{v_j: j\in [i-1]\}|\le r-1$ 
	for all $i\in [n]$. 
	It can be shown easily that $\chi_{DP}(G)\le col(G)$.
	Thus, for any graph $G$, 
	\equ{1-1}
	{
		\chi(G)\le \chi_\ell(G)\le \chi_{DP}(G)\le col(G)
	} 
	holds and the inequalities might be strict. 
For example, $\chi(K_{2,4})=2$ while $\chi_\ell(K_{2,4})=3$ in~\cite{Erdos1979}; $\chi(C_{2k})=\chi_\ell(C_{2k})=2$ while $\chi_{DP}(C_{2k})=3$ for all $k\ge 2$ and $k\in \N$ in~\cite{Dvorak2018}; and a result established in Section~\ref{sec3} that 
$\chi(C_4\vee K_1)=\chi_{DP}(C_4\vee K_1)=3$ while $col(C_4\vee K_1)=4$.

There are various studies on when the equalities in~(\ref{1-1}) will hold.
Thomassen \cite{Thom1994} proved $\chi_l(G)\le 5$ for 
every planar graph $G$ while 
Dvo\v{r}\'{a}k and Postle~\cite{Dvorak2018}
also observed that 
the DP-chromatic number of every planar graph is at most 5.
Another well-known result is due to Noel, Reed and Wu~\cite{Noel2015}.

\begin{theo}[\cite{Noel2015}]\label{Noel2015}
	Let $G$ be a graph with $n$ vertices. If $\chi(G) \ge (n-1)/2$, then $\chi_\ell(G) = \chi(G).$
\end{theo}

The condition $\chi(G) \ge (n-1)/2$ in Theorem~\ref{Noel2015} was shown to be the best possible due to many infinite classes of graphs. 
For example, for the complete $k$-partite graph $K_{5,2,\cdots,2}$,  
the condition in Theorem~\ref{Noel2015} fails and 
$
\chi_\ell(K_{5,2,\cdots,2})=k+1> \chi(K_{5,2,\cdots,2})
$
(see~\cite{Enomo2002} for more details).

Since any graph satisfying the requirement of Theorem~\ref{Noel2015} has a relatively large chromatic number, the join of any graph and a large enough complete graph will have its list chromatic number and chromatic number to be equal. 
Therefore, given any graph $G$, how large the order of the complete graph to join is enough becomes a natural question to be considered next.

To better describe this question, Enomoto, Ohba, Ota, and Sakamoto~\cite{Enomo2002} raised the following parameter:
$$Z_\ell(G)=\min\{m\in \mathbb{N}: \chi_\ell (G\vee K_m)=\chi(G\vee K_m)\}\text{ for any graph }G.$$
Note that for all integers $r\ge Z_\ell(G)$, $\chi_\ell (G\vee K_r)=\chi(G\vee K_r)$ holds. 

And Bernshteyn, Kostochka, and Zhu~\cite{Bern2017} further introduced three more parameters:
\begin{enumerate}
    \item 
$Z_{DP}(G)=\min\{m\in \mathbb{N}: \chi_{DP}(G\vee K_m)=\chi(G\vee K_m)\}$ for any graph $G$;
	\item
$Z_\ell(n)=\max_{|V(G)|=n}\{Z_\ell(G)\}$ for any $n\in \N$; and 
	\item 
$Z_{DP}(n)=\max_{|V(G)|=n}\{Z_{DP}(G)\}$ for any $n\in \N$.
\end{enumerate}

Similarly, for all integers $r\ge Z_{DP}(G)$, $\chi_{DP} (G\vee K_r)=\chi(G\vee K_r)$ holds. 

It is obvious that $Z_\ell(G)\le Z_{DP}(G)$ holds for any graph $G$ and $Z_\ell(n)\le Z_{DP}(n)$ holds for any $n\in \N$.
Furthermore, by Theorem~\ref{Noel2015}, $Z_\ell(n)$ has an upper bound $n-5$ for all $n \ge 5$ (mentioned in~\cite{Enomo2002}), while for DP coloring, Bernshteyn, Kostochka, and Zhu~\cite{Bern2017} has shown that $Z_{DP}(n)=\Omega(n^2)$ by the following two theorems.

\begin{theo}[\cite{Bern2017}]\label{Bern20172}
	Let $G$ be a graph with $n$ vertices and $m$ edges, then $Z_{DP}(G)\le 3m$. 
Hence $Z_{DP}(n)\le 3n^2/2$.
\end{theo}

\begin{theo}[\cite{Bern2017}]\label{Bern20173}
For any positive integer $n$, $Z_{DP}(n)\ge n^2/4-O(n)$. 
\end{theo}

In order to establish Theorem~\ref{Bern20172}, the authors of 
\cite{Bern2017} proved the following result.

\begin{theo}[\cite{Bern2017}]\label{Bern20171}
	Let $G$ be a $k$-colorable graph and $A$ be a complete graph with $|A|$ vertices.
If $\cov=(L,H)$ is a cover of $G\vee A$, where $|L(a)|\ge |A|+k$ for all $a\in V(A)$ and 
\equ{}
{
	|A|\ge 
\frac{3}{2}\sum_{v\in V(G)}\max\{d_G(v)-|L(v)|+|A|+1,0\},
}
then $G\vee A$ is $\cov$-colorable.
\end{theo}

In this paper, we shall improve the upper bounds of $Z_{DP}(G)$ 
to $\lceil 2.4m\rceil$ and $Z_{DP}(n)$ to $n^2-(n+3)/2$ by modifying the conditions in Theorem~\ref{Bern20171}.

\subsection{Main results  \label{sec1.2}}

Our main Theorem is as follows.

\begin{theo}\label{th1}
Let $G$ be a $k$-colorable graph and $A$ be a complete graph with $|A|$ vertices.
If $\cov=(L,H)$ is a cover of $G\vee A$, where $|L(a)|\ge |A|+\chi(G)$ for all $a\in V(A)$ and 
	\equ{th1-e0}
{
		|A|\ge 
		\var \sum_{v\in V(G)}\max\{d_G(v)-|L(v)|+|A|+k,0\},
}
then $G\vee A$ is $\cov$-colorable.
\end{theo}

We shall give a proof of Theorem~\ref{th1} in Section~\ref{sec2}, 
which is essentially a refinement and further derivation 
of the proof of Theorem~\ref{Bern20171} in~\cite{Bern2017}.

A graph $G$ is called {\it chordal} if for every cycle $C$ in $G$,
$G[V(C)]$ contains $3$-cycles. 
Due to Dirac~\cite{Dirac1961}, a graph $G$ 
is chordal if and only if  
there exists an ordering $v_1,v_2,\cdots,v_n$ of vertices in $G$,
called a {\it perfect elimination ordering},  
such that each vertex $v_i$ is simplicial in the subgraph of $G$
induced by $\{v_j: j\in [i]\}$.
Thus, for each chordal graph $G$, $col(G)=\chi(G)$ holds and all the equalities in (\ref{1-1}) naturally hold.

\begin{lemm}\label{le1-1}
Let $G$ be a graph with $n$ vertices.
If $G$ is chordal, then $\chi_{DP}(G)=\chi(G)$.
\end{lemm}

It can be verified that any graph $G$ with $n$ vertices and $\chi(G)\ge n-1$ has a clique with at least $n-1$ vertices, implying that $G$ is chordal. Then the next conclusion follows from Lemma~\ref{le1-1} directly. 

\begin{coro}\label{co1-1}
Let $G$ be a graph with $n$ vertices. If $\chi(G)\ge n-1$, then 
$\chi_{DP}(G)=\chi(G)$.
\end{coro}

When $\chi(G)=n-2$, 
we will prove in Section~\ref{sec3} that 
$\chi_{DP}(G)=\chi(G)$ holds
unless $G\cong C_4$.

\begin{theo}\label{th1-5}
(i) $Z_{DP}(C_4)=1$; 

(ii) For any graph $G$ with  $n$ vertices,  
if $\chi(G)=n-2$ and $G\not\cong C_4$, then $\chi_{DP}(G)=\chi(G)$.
\end{theo}

By Theorems~\ref{th1} and~\ref{th1-5}, we can obtain better upper bounds of both $Z_{DP}(G)$ and $Z_{DP}(n)$, the proofs of which are given in Section~\ref{sec3}.

\begin{theo}\label{th1-3}
For any graph $G$ with $n$ vertices, $m$ edges and chromatic number $k$,
\equ{th1-3-e1}
{
Z_{DP}(G)
\left \{
\begin{array}{ll} 
\le \left \lceil \frac{4(k+1)m}{2k+1} \right \rceil,
\quad &\mbox{if } 2\le k\le n-3;\\ 
=1, &\mbox{if }G\cong C_4;\\
=0, &\mbox{otherwise}.
\end{array}
\right.
}
\end{theo}

\begin{theo}\label{th1-4}
	For any integer $n\ge 2$,  $Z_{DP}(n)\le n^2-(n+3)/2$.
\end{theo}

\section{Proof of Theorem~\ref{th1}
\label{sec2}}

Given a cover $\cov=(L,H)$ of $G$, for any vertex-induced subgraph $K$ of $G$, let $\cov_K=(L_K, H_K)$ be a cover of $K$, 
where $L_K=L|_{V(K)}$ and $H_K=H[L(V(K))]$.
Thus, for each $u\in V(K)$, $L_K(u)=L(u)$.
For any $S\subseteq V(G)$, let $L(S)=\bigcup_{u\in S}L(u)$.

Now we are going to prove Theorem~\ref{th1}.

\noindent \textit{Proof of Theorem~\ref{th1}.}
Given a tuple $(k, G, A, \cov)$, where $G$ is $k$-colorable, $A$ is a complete graph with $|A|$ vertices and $\cov=(L,H)$ is a cover of graph $G\vee A$, for any $v\in V(G)$, let
$$
\sigma(v, k, G, A, \cov)=\max\{d_G(v)-|L(v)|+|A|+k,0\}
$$ 
and 
$$
\sigma(k, G, A, \cov)=\sum\limits_{v\in V(G)}\sigma(v, k, G, A, \cov).
$$

Assume that there is a counterexample of the statement, the tuple $(k, G, A, \cov)$, which minimizes $k$, then $|V(G)|$, then $|A|$. Obviously, $k\ge \chi(G)$.

\Clm{cl1}
{$|V(G)|\ge 2$.}

\proof
If $|V(G)|=0$, then $|L(a)|\ge |A|>d_{G\vee A}(a)$ for each $a\in A$, 
implying that $G\vee A$ is $\cov$-colorable.

Suppose that $|V(G)|=1$. Then $k\ge \chi(G)=1$.
The condition of (\ref{th1-e0}) implies that $|L(v)|\ge k\ge 1$
	for $v\in V(G)$.
As $|V(G)|=1$, $d_{G\vee A}(a)=|A|<|L(a)|$ for each $a\in V(A)$, 
implying that $G\vee A$ is $\cov$-colorable, a contradiction.

Thus Claim~\ref{cl1} holds.
\claimend

\Clm{cl2}
{$|A|\ge 1$.}

\proof
Assume $|A|=0$. Then for all $v\in V(G)$, $|L(v)|\ge d_G(v)+k\ge d_G(v)+1$ holds, 
implying that $G\vee A$ is $\cov$-colorable, a contradiction. 
Hence Claim~\ref{cl2} holds.
\claimend

\Clm{cl3}
{
	$G_v\vee A$ is $\cov_{G_v\vee A}$-colorable for all $v\in V(G)$.
}

\proof
For any $v\in V(G)$, $k\ge \chi(G)\ge \chi(G_v)$ holds. So we consider the tuple $(k, G_v, A, \cov_{G_v\vee A})$, where $L_{G_v\vee A}(u)=L(u)$ for all $u\in V(G_v\vee A)$.
Note that for each $a\in V(A)$,
$|L_{G_v\vee A}(a)|=|L(a)|\ge |A|+\chi(G)\ge |A|+\chi(G_v)$,
and for each $u\in V(G_v)$, 
\aln{th1-e1}
{
	\sigma(u, k, G_v, A, \cov_{G_v\vee A})&=\max\{d_{G_v}(u)-|L_{G_v\vee A}(u)|+|A|+k,0\}
	\nonumber\\
	&\le \max\{d_G(u)-|L(u)|+|A|+k,0\}
	\nonumber\\
	&= \sigma(u, k, G, A, \cov).
}
Thus
\equ{th1-e2}
{
	|A|\ge \var\sigma(k, G, A, \cov)
	\ge \var\sum_{u\in V(G)\setminus \{v\}}\sigma(u,k, G, A, \cov)
	\ge \var\sigma(k, G_v, A, \cov_{G_v\vee A}).
}

By the assumption on the minimality of $|V(G)|$, $G_v\vee A$ is $\cov_{G_v\vee A}$-colorable. 
Thus Claim~\ref{cl3} holds.
\claimend

\Clm{cl4}
{
	For each $v\in V(G)$,  $\sigma(v, k, G, A, \cov)=d_G(v)-|L(v)|+|A|+k\ge k$.
}

\proof
Assume there exists $v_0 \in V(G)$ with $\sigma(v_0, k, G, A, \cov)\le k-1$. Then 
$$
|L(v_0)|\ge d_G(v_0)+|A|+1=d_{G\vee A}(v_0)+1.
$$
Since $G_{v_0}\vee A$ is $\cov_{G_{v_0}\vee A}$-colorable by Claim 3, 
$|L(v_0)|\ge d_{G\vee A}(v_0)+1$ implies that 
$G\vee A$ is $\cov$-colorable, a contradiction. 
Hence Claim~\ref{cl4} holds.
\claimend

\Clm{cl5}
{
	There does not exist a cycle $x_1-y_1-\cdots-x_\ell-y_\ell-x_1$ in $H$, where 
	\begin{enumerate}
		\item $2\le \ell\le 2k+1$,
		\item $x_i\in L(V(A))$ and $y_i\in L(V(G))$ for all $1\le i\le \ell$, and 
		\item $\{x_1,\cdots,x_\ell\}$ is an independent set in $H$.
	\end{enumerate}
}

\tikzstyle{cblue}=[circle, draw, thin,fill=blue!20, scale=0.5]
\begin{figure}[!ht]
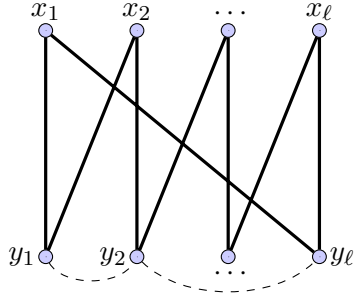

	\tikzp{1.5}
	{
		\foreach \place/\y in {{(-1.2,1)/1}, {(-0.4,1)/2},{(0.4,1)/3}, {(1.2,1)/4}, {(-1.2,-1)/5}, {(-0.4,-1)/6},{(0.4,-1)/7},{(1.2,-1)/8}}   
		\node[cblue] (b\y) at \place {};
		
		\filldraw[black] (b1) circle (0pt)node[anchor=south] {$x_1$};
		\filldraw[black] (b2) circle (0pt)node[anchor=south] {$x_2$};
		\filldraw[black] (b3) circle (0pt)node[anchor=south] {~$\cdots$};
		\filldraw[black] (b4) circle (0pt)node[anchor=south] {$x_\ell$};
		\filldraw[black] (b5) circle (0pt)node[anchor=east] {$y_1$};
		\filldraw[black] (b6) circle (0pt)node[anchor=east] {$y_2$};
		\filldraw[black] (b7) circle (0pt)node[anchor=north] {~$\cdots$};
		\filldraw[black] (b8) circle (0pt)node[anchor=west] {$y_\ell$};		
		
		\draw[black, very thick] (b1) -- (b5) -- (b2)--(b6) -- (b3) -- (b7)--(b4)--(b8) --(b1);  
		\draw [black, dashed, bend right=60, looseness=0.85] (b5) to (b6);
		\draw [black, dashed, bend right=45, looseness=0.85] (b6) to (b8);
	}
	\caption{An example of Claim 5}
	\label{f1}
\end{figure}

\proof
Assume there exists a cycle $x_1-y_1-\cdots-x_\ell-y_\ell-x_1$ in $H$ satisfying (i)-(iii). Note that all the vertices on the cycle are pairwise distinct.

Let $X=\{x_i:1\le i\le \ell\}$,
$x_i\in L(a_i)$ for all $i:1\le i\le \ell$ and $A_0=\{a_i:1\le i\le \ell\}$.
As $X$ is an independent set of $H$, 
we have $|A_0|=\ell$.

Let $A'=A\setminus A_0$ and $\cov'=(L',H')$ be a cover of $G\vee A'$,
where $L'(u)=L(u)\setminus N_H(X)$ for each $u\in V(G\vee A')$
and 
$H'$ is the subgraph of $H$ induced by $\bigcup_{u\in V(G\vee A')}L'(u)$.

Consider the tuple $(k, G, A', \cov')$.
For each $u\in V(G\vee A')$, 
$|L'(u)|\ge |L(u)|-\ell$ holds.
Thus, for each $a\in V(A')$, 
\equ{th1-e3}
{
		|L'(a)|\ge |L(a)|-\ell
		\ge |A|+\chi(G)-\ell= |A'|+\chi(G).
} 

Let $Y=\{y_i: 1\le i\le \ell\}$.
We are now going to show that
for any vertex $v$ in $V(G)$,
$|L'(v)|\ge |L(v)|-\ell+|L(v)\cap Y|$.
Let  $q=|L(v)\cap Y|$.
Assume that $L(v)\cap Y=\{y_{t_1},y_{t_2},\cdots,y_{t_q}\}$,
where $1\le t_1<t_2<\cdots<t_q\le \ell$. 
For each $i$ with $1\le i\le \ell$, 
as $\{y_i,y_{i+1}\}\subseteq N_H(x_{i+1})$
and $E_H(L(v),L(a_{i+1}))$ is a matching of $H$,
where $x_{\ell+1}=x_1$, $y_{\ell+1}=y_1$ and $a_{\ell+1}=a_1$, 
we have $|L(v)\cap\{y_i,y_{i+1}\}|\le 1$.
Thus, 
$0\le q\le \lfloor\ell/2\rfloor$ and every two numbers in $t_1,\cdots, t_q$ 
are not consecutive,  and in particular, $t_1\ne 1$ whenever $t_q=\ell$.
It implies that 
\equ{th1-e4}
{
	|X\cap N_H(L(v)\cap Y)|\ge |X'|
	=2q,
}
where $X'=\{x_{t_j},x_{t_j+1}: j=1,2,\cdots,q \}$.
Then
\eqn{th1-e5}
{
	|L(v)\cap N_H(X)|
	&\le &|L(v)\cap N_H(X')|
	+|L(v)\cap N_H(X\setminus X')|
	\nonumber \\
	&\le &|L(v)\cap Y|+|X\setminus X'|
	\nonumber \\
	&=&q+(\ell-2q)=\ell-q,
}
by which $|L'(v)|=|L(v)|-|L(v)\cap N_H(X)|\ge |L(v)|-(\ell-q)$.
Thus, $|A'|-|L'(v)|\le |A|-|L(v)|-q$ and 
\eqn{th1-e6}
{
	\sigma(v, k, G, A', \cov')
	&=&\max\{d_G(v)-|L'(v)|+|A'|+k,0\}
	\nonumber\\
	&\le &\max\{d_G(v)-|L(v)|+|A|+k-q,0\}
	\nonumber\\
	&=& d_G(v)-|L(v)|+|A|+k-q
	\nonumber\\
	&=&\sigma(v, k, G, A, \cov)-|L(v)\cap Y|,
}
where the second last equality holds due to Claim 4 and $0\le q\le \lfloor \ell/2\rfloor \le k$.
Hence 
\eqn{th1-e7}
{
	\sigma(k, G, A', \cov')
	&=&\sum_{v\in V(G)}\sigma(v, k, G, A', \cov')
	\nonumber \\ 
	&\le &\sum_{v\in V(G)}\sigma(v, k, G, A, \cov)-\sum_{v\in V(G)} |L(v)\cap Y|
	\nonumber \\ 
	&=&\sigma(k, G, A, \cov)-\ell
}
and
\eqn{th1-e8}
{
	|A'|&=&|A|-\ell
	\nonumber\\
	&\ge& \var\sigma(k, G, A, \cov)-\ell
	\nonumber\\
	&\ge& \var(\sigma(k, G, A', \cov')+\ell)-\ell
	\nonumber\\
	&>& \var\sigma(k, G, A', \cov').
}

By the assumption on the minimality of $|A|$, $G\vee A'$ is $\cov'$-colorable, whereas the union of $X$ and any $\cov'$-coloring of $G\vee A'$ is an $\cov$-coloring of $G\vee A$, a contradiction. Thus Claim~\ref{cl5} holds.
\claimend

\Clm{cl6}
{ 
	For any $a\in V(A)$, $x\in L(a)$ and $v\in V(G)$, $N_H(x)\cap L(v)\neq \emptyset.$
}

\proof
Assume there exist $a_0\in V(A)$, $x_0\in L(a_0)$, and $v_0\in V(G)$, such that $N_H(x_0)\cap L(v_0)=\emptyset.$

Let $A'=A\setminus \{a_0\}$ and $\cov'=(L', H')$ be a cover of $G\vee A'$, where
$L'(u)=L(u)\setminus N_H(x_0)$ for each $u\in V(G\vee A')$
and $H'$ is the subgraph of $H$ induced by 
$\bigcup_{u\in V(G\vee A')}L'(u)$.
Consider the tuple $(k, G, A', \cov')$.

Then $|L'(u)|\ge |L(u)|-1$ holds for each $u\in V(G\vee A')$, 
and especially, $|L'(v_0)|= |L(v_0)|$ holds.
Consequently, $|L'(a)|\ge |L(a)|-1\ge |A'|+\chi(G)$ holds for each $a\in V(A')$.
For each $v\in V(G)$, 
\eqn{th1-e9}
{
	\sigma(v, k, G, A', \cov')
	&=&\max\{d_G(v)-|L'(v)|+|A'|+k,0\}
	\nonumber \\
	&\le & \max\{d_G(v)-|L(v)|+1+|A|-1+k,0\}
	\nonumber \\
	&=&\sigma(v, k, G, A, \cov),
}
and especially, $\sigma(v_0, k, G, A', \cov')=\sigma(v_0, k, G, A, \cov)-1$,
implying that 
$$
\sigma(k, G, A', \cov') \le \sigma(k, G, A, \cov)-1
$$ 
and 
\eqn{th1-e10}
{
	|A'|&=&|A|-1
	\nonumber\\
	&\ge& \var\sigma(k, G, A, \cov)-1
	\nonumber\\
	&\ge& \var(\sigma(k, G, A', \cov')+1)-1
	\nonumber\\
	&>& \var\sigma(k, G, A', \cov').
}

By the assumption on the minimality of $|A|$, $G\vee A'$ is $\cov'$-colorable, 
whereas the union of $\{x_0\}$ and any $\cov'$-coloring of $G\vee A'$ is an $\cov$-coloring of $G\vee A$, a contradiction. Thus Claim~\ref{cl6} holds.
\claimend

\Clm{cl7}
{  
	For each $v\in V(G)$, $|L(v)| \ge |A|+\chi(G)$ and 
$d_G(v)\ge \chi(G)\ge 2$.
}

\proof
By Claim 6, for any $v\in V(G)$, $a\in V(A)$ and $x\in L(a)$, 
$N_H(x)\cap L(v)\neq \emptyset$. Since $E_H(L(a), L(v))$ is a matching, we have $|L(v)|\ge |L(a)|\ge |A|+\chi(G)$.
Thus, by Claim 4,
$$d_G(v)+|A|+k\ge|L(v)|+k\ge |A|+\chi(G)+k.$$
This implies $d_G(v)\ge \chi(G)\ge 1$ for all $v\in V(G)$, 
thus $d_G(v)\ge \chi(G)\ge 2$ further holds.

Hence  Claim~\ref{cl7} holds.
\claimend

\Clm{cl08}
{
	Any proper $\chi(G)$-coloring of $G$ has 
		two color classes of sizes at least $2$.
}

\proof 
Suppose Claim~\ref{cl08} fails.
Then there exists a proper 
$r$-coloring of $G$, where $r=\chi(G)$,
such that its color classes are $U_1,U_2,\cdots,U_r$
with $|U_2|= \cdots= |U_r|=1$.

Then $G[U_2\cup\cdots \cup U_r]$ is isomorphic to $K_{r-1}$.
Otherwise, $\chi(G)<r$, a contradiction.
But, as $G[U_2\cup\cdots \cup U_r]\cong K_{r-1}$,
each vertex in $U_1$ is a simplicial vertex,
implying that $G$ has a perfect elimination ordering
and $G$ is chordal. 

By Claim~\ref{cl7} and the given condition, 
we have $|L(u)|\ge |A|+\chi(G)$ for each $u\in V(G\vee A)$.

Since $G\vee A$ is also chordal, $\chi_{DP}(G\vee A)=\chi(G\vee A)=\chi(G)+|A|$ by Lemma~\ref{le1-1},
and therefore $G\vee A$ has an $\cov$-coloring,
a contradiction.
	
Thus Claim~\ref{cl08} holds.
\claimend

\Clm{cl10}
{
For any vertex $v\in V(G)$ and any independent set 
$I\subseteq L(V(A))$ in $H$,
if $|I|\le k+d_G(v)$, then 
$|I|+|L(v)| \le |L(v)\cap N_H(I)|+d_G(v)+|A|$.
}

\proof Assume that $v\in V(G)$ and 
$I\subseteq L(V(A))$ is an independent set in $H$ 
with $|I|\le k+d_G(v)$.
Due to Claim~\ref{cl4}, the claim holds when $I=\emptyset$.
Now assume that $I\ne \emptyset$.

Let $A_0=\{a\in A:  L(a)\cap I\ne\emptyset\}$.
As $I$ is an independent set, $|A_0|=|I|>0$.
Recall that $G_v$ is the subgraph $G\setminus \{v\}$.
Let $A'=A\setminus A_0$ and $\cov'=(L',H')$ be the cover of $G_v\vee A'$, where
$L'(u)=L(u)\setminus N_H(I)$ for all $u\in V(G_v\vee A')$
and $H'$ is the subgraph of $H$ induced by 
$\bigcup_{u\in V(G_v\vee A')}L'(u)$.

Consider the tuple $(k,G_v,A',\cov')$ as $k\ge \chi(G)\ge \chi(G_v)$.
For each $a\in A'$, $|L'(a)|\ge 
|L(a)|-|I|\ge |A|+\chi(G)-|I|\ge |A'|+\chi(G_v)$.
For each $u\in V(G_v)$, $|L'(u)|\ge |L(u)|-|I|$ and 
\eqn{cl10-e1}
{
\sigma(u, k,G_v,A',\cov')
&=&\max\{d_{G_v}(u)-|L'(u)|+|A'|+k,0\}
\nonumber \\
&\le &\max\{d_{G}(u)-|N_G(v)\cap \{u\}|-(|L(u)|-|I|)+|A|-|I|+k,0\}
\nonumber \\
&=&\max\{d_{G}(u)-|N_G(v)\cap \{u\}|-|L(u)|+|A|+k,0\}
\nonumber \\
&=&\sigma(u, k,G,A,\cov)-|N_G(v)\cap \{u\}|,
}
where the last equality follows from the fact $\sigma(u, k,G,A,\cov)\ge k\ge 2$ 
due to Claims~\ref{cl4} and~\ref{cl7}.
Thus, 
\eqn{cl10-e2}
{
\sigma(k,G_v,A',\cov')&\le &\sum_{u\in V(G_v)} 
\left ( \sigma(u, k,G,A,\cov)-|N_G(v)\cap \{u\}|\right )
\nonumber \\ 
&=&\sum_{u\in V(G_v)}\sigma(u, k,G,A,\cov)
-d_G(v)
\nonumber \\ 
&=&\sigma(k,G,A,\cov)-\sigma(v, k,G,A,\cov)
-d_G(v)
\nonumber \\ 
&\le &\sigma(k,G,A,\cov)-k-d_G(v),
}
where the last inequality follows form the fact $\sigma(v, k,G,A,\cov)\ge k$ 
by Claim~\ref{cl4}, 
and
\eqn{cl10-e3}
{
|A'|&=&|A|-|I|
\nonumber \\ 
&\ge &\var \sigma(k,G,A,\cov) -|I|
\nonumber \\ 
&\ge &\var \left (\sigma(k,G_v,A',\cov')+ k+d_G(v) \right ) -|I|
\nonumber \\ 
&\ge &\var \sigma(k,G_v,A',\cov'),
}
where the last inequality follows from the condition 
that $|I|\le k+d_G(v)$.

By the assumption on the minimality of $|V(G)|$, 
$G_v\vee A'$ has an $\cov'$-coloring, i.e., an independent set 
$Q_1$ of $H'$ with $|Q_1|=|V(G)|+|A|-|I|-1$.
Clearly, $Q_2:=Q_1\cup I$ is an independent set of $H$ 
with $|Q_2|=|V(G)|+|A|-1$ and $Q_2\cap L(v)=\emptyset$.

As $G\vee A$ has no $\cov$-coloring by assumption,  
we have $L(v)\subseteq N_H(Q_2)$,
 implying that 
\eqn{cl10-e4}
{
|L(v)| &=&|L(v)\cap N_H(Q_1\cup I)| 
\nonumber \\
&\le &|L(v)\cap N_H(I)|+|L(v)\cap N_H(Q_1\cap L(V(A')))|+
|L(v)\cap N_H(Q_1\cap L(V(G_v)))|
\nonumber \\
&\le &|L(v)\cap N_H(I)|+|A'|+d_G(v)
\nonumber \\
&= &|L(v)\cap N_H(I)|+|A|+d_G(v)-|I|,
}
by which Claim~\ref{cl10} follows.
\claimend


\Clm{cl11}
{
For any $v\in V(G)$, if $\sigma(v,k,G,A,\cov)=k$,
then $N_H(y)\cap L(V(A))$ is a clique of $H$ for each $y\in L(v)$.
}

\proof 
Assume that $\sigma(v,k,G,A,\cov)=k$ for $v\in V(G)$. 
Then $|L(v)|=d_G(v)+|A|$.
Suppose that $N_H(y)\cap L(V(A))$ is not a clique of $H$ for some $y\in L(v)$.
Then there exists an independent set $I=\{z_1,z_2\}\subseteq N_H(y)\cap L(V(A))$. 
Clearly, $|I|=2\le k+d_G(v)$  by Claim~\ref{cl7}.
Thus, due to Claim~\ref{cl10},  $|N_H(I)\cap L(v)|\ge |I|=2$. 
But, $I\subseteq N_H(y)$ implies that $N_H(I)\cap L(v)=\{y\}$,
a contradiction.

Thus Claim~\ref{cl11} holds. 
\claimend

Claim~\ref{cl11} will be applied in the proof of Claim~\ref{cl12}.

Due to Claim~\ref{cl08},  
there are two disjoint independent sets $U_0$ and $U_1$ in $G$, 
such that $|U_i|\ge 2$
and $\chi(G\setminus U_i)=\chi(G)-1$ for both $i=0,1$.

\Clm{cl9-0}
{
For any path $P: x_1-y_1-x_2-y_2-\cdots - y_i-x_{i+1}$  in $H$, 
where $0\le i\le 2k$, 
if $y_j\in L(U_{j(\text{mod 2})})$ 
for each $j\in [i]$ and
$\{x_1,x_2,\cdots,x_{i+1}\}\subseteq L(V(A))$ 
	is an independent set,
then $\{x_1,x_2,\cdots,x_{i}\}\cup 
	(N_H(x_{i+1})\cap L(U_{(i+1)(\text{mod 2})}))$
	is independent in $H$ with size $i+|U_{(i+1)(\text{mod 2})}|$.
}

\proof Let $Y=\{y_j: 1\le j\le i\}$.

We first show that 
$N_H(x_{i+1})\cap L(U_{(i+1)(\text{mod 2})})\cap Y=\emptyset$.
Otherwise, there exists $y_t\in Y$ such that 
$y_t\in N_H(x_{i+1})\cap L(U_{(i+1)(\text{mod 2})})
$, 
where $1\le t\le i$.
Clearly, $t<i$,  as $y_i\in L(U_{i(\text{mod 2})})$. Therefore, $x_{t+1}-y_{t+1}-\cdots-x_i-y_i-x_{i+1}-y_t-x_{t+1}$ is a cycle in $H$,
as shown in Figure~\ref{f3}, 
a contradiction to Claim~\ref{cl5}. 
Thus, the conclusion holds.

\tikzstyle{cblue}=[circle, draw, thin,fill=blue!20, scale=0.5]
\begin{figure}[!ht]
	\tikzp{1.5}
	{
		\foreach \place/\y in {{(-0.4,1)/2},{(0.4,1)/3}, {(1.2,1)/4}, {(-1.2,-1)/5}, {(-0.4,-1)/6},{(0.4,-1)/7}}   
		\node[cblue] (b\y) at \place {};
		
		\filldraw[black] (b2) circle (0pt)node[anchor=south] {$x_{t+1}$};
		\filldraw[black] (b3) circle (0pt)node[anchor=south] {~$\cdots$};
		\filldraw[black] (b4) circle (0pt)node[anchor=south] {$x_{i+1}$};
		\filldraw[black] (b5) circle (0pt)node[anchor=north] {$y_t$};
		\filldraw[black] (b6) circle (0pt)node[anchor=north] {~$\cdots$};
		\filldraw[black] (b7) circle (0pt)node[anchor=north] {$y_i$};	
		
		\draw[black, very thick] (b5) -- (b2)--(b6) -- (b3) -- (b7)--(b4)--(b5);  
	}
	\caption{A cycle}
	\label{f3}
\end{figure}

Let $Q=\{x_1,x_2,\cdots,x_{i}\}\cup 
	(N_H(x_{i+1})\cap L(U_{(i+1)(\text{mod 2})}))$.
Clearly, if $i=0$, then 
$Q=N_H(x_{i+1})\cap L(U_{(i+1)(\text{mod 2})})$.

As $U_{(i+1)(\text{mod 2})}$ is an independent set, 
$N_H(x_{i+1})\cap L(U_{(i+1)(\text{mod 2})})$ is
an independent set in $H$.
By Claim~\ref{cl6},
$|N_H(x_{i+1})\cap L(u)|=1$
for each $u\in  U_{(i+1)(\text{mod 2})}$,
implying that $|Q|=i+|U_{(i+1)(\text{mod 2})}|$.
Clearly, Claim~\ref{cl9-0} holds when $i=0$.

Now suppose that $i\ge 1$ and that $Q$ is not independent.
Since both $\{x_1,x_2,\cdots,x_{i}\}$ and 
$N_H(x_{i+1})\cap L(U_{(i+1)(\text{mod 2})})$
are independent, 
there exists
$y_{i+1}\in N_H(x_t)\cap (N_H(x_{i+1})\cap 
	L(U_{(i+1)(\text{mod 2})}))$ 
for some $1\le t\le i$.
Since $N_H(x_{i+1})\cap L(U_{(i+1)(\text{mod 2})})\cap Y=\emptyset$,
$y_{i+1}\notin Y$.
Then $x_t-y_t-\cdots-x_{i+1}-y_{i+1}-x_t$ is a cycle in $H$,
as shown in Figure~\ref{f4}, 
a contradiction to Claim~\ref{cl5}.

Thus, Claim~\ref{cl9-0} holds.
\claimend

\tikzstyle{cblue}=[circle, draw, thin,fill=blue!20, scale=0.5]
\begin{figure}[!ht]
	\tikzp{1.5}
	{
		\foreach \place/\y in {{(-0.4,1)/2},{(0.4,1)/3}, {(1.2,1)/4}, {(-0.4,-1)/6},{(0.4,-1)/7},{(1.2,-1)/8}}   
		\node[cblue] (b\y) at \place {};
		
		\filldraw[black] (b2) circle (0pt)node[anchor=south] {$x_t$};
		\filldraw[black] (b3) circle (0pt)node[anchor=south] {~$\cdots$};
		\filldraw[black] (b4) circle (0pt)node[anchor=south] {$x_{i+1}$};
		\filldraw[black] (b6) circle (0pt)node[anchor=east] {$y_t$};
		\filldraw[black] (b7) circle (0pt)node[anchor=north] {~$\cdots$};
		\filldraw[black] (b8) circle (0pt)node[anchor=west] {$y_{i+1}$};		
		
		\draw[black, very thick] (b2)--(b6) -- (b3) -- (b7)--(b4)--(b8)--(b2);  
	}
	\caption{A cycle}
	\label{f4}
\end{figure}

\Clm{cl9}
{There exists a path $x_1-y_1-x_2-y_2-\cdots - y_{2k+1}-x_{2k+2}$ in $H$
	such that $y_j\in L(U_{j(\text{mod 2})})$ for each $j\in [2k+1]$
	and
	$\{x_1,x_2,\cdots,x_{2k+2}\}\subseteq L(V(A))$ 
	is an independent set.
}

\proof
	Suppose Claim~\ref{cl9} fails. Let 
	$P: x_1-y_1-x_2-y_2-\cdots - y_i-x_{i+1}$ be a longest path in $H$, 
where $0\le i\le 2k$, 
such that $y_j\in L(U_{j(\text{mod 2})})$ 
for each $j\in [i]$ and
$\{x_1,x_2,\cdots,x_{i+1}\}\subseteq L(V(A))$ 
	is an independent set.
	Such a path exists, at least for $i=0$.
Note that all the vertices on $P$ are pairwise distinct.

Let $x_j\in L(a_j)$ for each $j\in [i+1]$, where $a_j\in V(A)$.
As $\{x_1,x_2,\cdots,x_{i+1}\}$ is an independent set in $H$,
$a_1,a_2,\cdots,a_{i+1}$ are pairwise distinct.

Due to the assumption of $P$,
there are no $x_{i+2}\in L(V(A)) 
\setminus \{x_1,x_2,\cdots,x_{i+1}\}$
and $y_{i+1}\in N_H(x_{i+1})\cap N_H(x_{i+2})\cap L(U_{(i+1)(\text{mod 2})})$,
such that $\{x_1,x_2,\cdots,x_{i+2}\}$ is independent in $H$.
Thus by Claim~\ref{cl9-0},
 for each $x\in L(V(A)) 
\setminus \{x_1,x_2,\cdots,x_{i+1}\}$, 
either $N_H(x)\cap \{x_1,x_2,\cdots,x_{i+1}\}\ne \emptyset$ or
$N_H(x)\cap N_H(x_{i+1})\cap L(U_{(i+1)(\text{mod 2})})=\emptyset$.

It implies that for any $a\in V(A)\setminus \{a_1,a_2,\cdots,a_{i+1}\}$,
\equ{cl9-e1}
{
	L(a)\cap 
	N_H(N_H(x_{i+1})\cap 
	L(U_{(i+1)(\text{mod 2})}))\subseteq L(a)\cap N_H(x_1,x_2,\cdots,x_{i+1}).
}

Let $Q=\{x_1,x_2,\cdots,x_{i}\}\cup 
	(N_H(x_{i+1})\cap L(U_{(i+1)(\text{mod 2})}))$.
	Clearly, if $i=0$, then $Q=N_H(x_{i+1})\cap L(U_{(i+1)(\text{mod 2})})$. 
Let $G'=G\setminus U_{(i+1)\text{(mod 2)}}$, 
$A'=A\setminus \{a_1,\cdots,a_{i}\}$, 
and $\cov'=(L', H')$ be the cover of $G'\vee A'$, where
$L'(u)=L(u)\setminus N_H(Q)$ for each $u\in V(G'\vee A')$ 
and $H'$ is the subgraph of $H$ induced by $\bigcup_{u\in V(G'\vee A')}L'(u)$. Thus, $A'=A$ if $i=0$.

We are now going to show that 
$G'\vee A'$ is $\cov'$-colorable. 
Consider the tuple 
$(k, G', A', \cov')$.
Note that $L(a_{i+1})\cap N_H(N_H(x_{i+1})\cap L(U_{(i+1)(\text{mod 2})}))
=\{x_{i+1}\}$, implying that 
$$
|L(a_{i+1})\cap N_H(Q)|   
=|L(a_{i+1})\cap N_H(x_1,x_2,\cdots,x_{i})|+1
\le i+1.
$$
For any $a\in V(A)\setminus \{a_1,a_2,\cdots,a_{i+1}\}$, 
by (\ref{cl9-e1}),
$$
|L(a)\cap N_H(Q)| 
\le |L(a)\cap N_H(x_1,x_2,\cdots,x_{i+1})|
\le i+1. 
$$
Hence, for each $a\in V(A')$, 
$
|L'(a)|\ge |L(a)|-i-1\ge |A|+\chi(G)-i-1=|A'|+\chi(G'),
$
where the last equality follows from the fact $\chi(G')=\chi(G)-1$.

For each $v\in V(G')$, 
$$
|L'(v)|=|L(v)\setminus N_H(Q)|\ge |L(v)|
-i-|N_G(v)\cap U_{(i+1)\text{(mod 2)}}|,
$$
implying that 
\eqn{th1-eq14}
{
	d_{G'}(v)-|L'(v)|
	&\le& d_{G}(v)-|N_G(v)\cap U_{(i+1)\text{(mod 2)}}|
	-|L(v)|+i+|N_G(v)\cap U_{(i+1)\text{(mod 2)}}|
	\nonumber \\
	&=&d_{G}(v)-|L(v)|+i
}
and
\eqn{th1-eq15}
{
	\sigma(v, k, G', A', \cov')
	&=&\max\{d_{G'}(v)-|L'(v)|+|A'|+k,0\}
	\nonumber \\
	&\le &\max\{d_{G}(v)-|L(v)|+i+|A|-i+k,0\}
	\nonumber \\
	&=&\max\{d_{G}(v)-|L(v)|+|A|+k,0\}
	\nonumber \\
	&=&\sigma(v, k, G, A, \cov).
}

Since $|U_{(i+1)\text{(mod 2)}}|\ge 2$,
by Claim~\ref{cl4} and (\ref{th1-eq15}), 
\eqn{cl9-e2}
{
	\sigma(k, G, A, \cov)
	&=&
	\sum_{v\in V(G')}\sigma(v,k, G, A, \cov)
	+\sum_{v\in U_{(i+1)\text{(mod 2)}}}\sigma(v,k, G, A, \cov)
	\nonumber \\ 
	&\ge & 
	\sigma(k, G', A', \cov')+2k.
}
Since $i\le 2k$,
\eqn{th1-eq16}
{
	|A'|&=&|A|-i
	\nonumber\\
	&\ge &\var\sigma(k, G, A, \cov)-i
	\nonumber\\
	&\ge & \var(\sigma(k, G', A', \cov')+2k)-2k
	\nonumber \\
	&> &\var\sigma(k, G', A', \cov'). 
}
By the assumption on the minimality of $|V(G)|$, 
$G'\vee A'$ is $\cov'$-colorable. 

Due to Claim~\ref{cl9-0}, $Q$ is independent in $H$. 
Since $|Q\cap L(u)|=1$ for each 
$u\in (A\setminus A')\cup U_{(i+1)\text{(mod 2)}}$,
the definition of $\cov'$ implies that 
the union of $Q$ 
and any $\cov'$-coloring of $G'\vee A'$ is an $\cov$-coloring of $G\vee A$, a contradiction.

Thus Claim~\ref{cl9} holds.
\claimend

\begin{figure}[!ht]
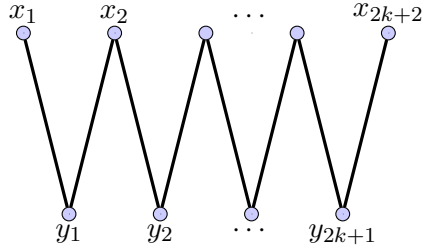

	\tikzp{1.5}
	{
		\foreach \place/\y in {{(-1.2,0.8)/1}, {(-0.4,0.8)/2},{(0.4,0.8)/3},
			{(1.2,0.8)/4}, {(-0.8,-0.8)/5}, {(0,-0.8)/6},{(0.8,-0.8)/7},{(2,0.8)/8}, {(1.6,-0.8)/9}}   
		\node[cblue] (b\y) at \place {};
		
		\filldraw[black] (b1) circle (0pt)node[anchor=south] {$x_1$};
		\filldraw[black] (b2) circle (0pt)node[anchor=south] {$x_2$};
		\filldraw[black] (0.8,0.8) circle (0pt)node[anchor=south] {$\cdots$};
		\filldraw[black] (b5) circle (0pt)node[anchor=north] {$y_1$};
		\filldraw[black] (b6) circle (0pt)node[anchor=north] {$y_2$};
		\filldraw[black] (b7) circle (0pt)node[anchor=north] {$\cdots$};
		\filldraw[black] (b8) circle (0pt)node[anchor=south] {$x_{2k+2}$};
		\filldraw[black] (b9) circle (0pt)node[anchor=north] {$y_{2k+1}$};		
		
		\draw[black, very thick] (b1) -- (b5) -- (b2)--(b6) -- (b3) -- (b7) -- (b4)--(b9)--(b8);  
	}
	\caption{A path $x_1-y_1-x_2-\cdots -y_{2k+1}-x_{2k+2}$ in $H$, where $\{x_1,\cdots,x_{2k+2}\}\subseteq L(V(A))$ is an independent set
		and $y_j\in L(U_{j\text{(mod 2)}})$}
	\label{f2}
\end{figure}

By Claim~\ref{cl9}, there exists a path 
$x_1-y_1-x_2-\cdots -y_{\vard}-x_{\varn}$
in $H$, as shown in Figure~\ref{f2},
where $\{x_1,\cdots,x_{\varn}\}\subseteq L(V(A))$ is an independent set
in $H$
and $y_j\in L(U_{j\text{(mod 2)}})$ for each $j=1,2,\cdots,\vard$. 

Let $X=\{x_j:j\in [2k+2]\}$ and $Y=\{y_j:j\in [2k+1]\}$.

\Clm{cl12}
{
For any $v\in V(G)$, $|L(v)\cap Y|\le k+1$, $|L(v)\cap N_H(X)|\le 2k+2-|L(v)\cap Y|$
and $\sigma(v,k,G,A,\cov)-|L(v)\cap Y|\ge 0$.
}

\proof For any $j:1\le j\le 2k+1$, $y_j\in L(U_{j\text{(mod 2)}})$, implying that 
for $1\le j\le 2k$, 
$|L(v)\cap \{y_j,y_{j+1}\}|=0$ if $v\notin U_0\cup U_1$,
and $|L(v)\cap \{y_j,y_{j+1}\}|\le 1$ otherwise. 
Thus, $|L(v)\cap Y|\le \lceil |Y|/2\rceil =\lceil (2k+1)/2\rceil =k+1$.

Note that if $L(v)\cap Y=\{y_{r_1}, \cdots,y_{r_t}\}$, then 
$X_0:=\{x_{r_j},x_{r_{j+1}}: j\in [t]\}\subseteq X\cap N_H(L(v)\cap Y)$ and 
\eqn{cl12-e1}
{
|L(v)\cap N_H(X)|&\le & |L(v)\cap Y|+|L(v)\cap N_H(X\setminus X_0)|
\nonumber \\
&\le & |L(v)\cap Y|+|X\setminus X_0|
\nonumber \\
&= & |L(v)\cap Y|+|X|-2|L(v)\cap Y|
\nonumber \\
&= &2k+2-|L(v)\cap Y|.
}

Now we are going to show that  if $L(v)\cap Y\ne \emptyset$, 
then $\sigma(v,k,G,A,\cov)\ge k+1$.
Otherwise, $\sigma(v,k,G,A,\cov)=k$ due to Claim~\ref{cl4}.
As $L(v)\cap Y\ne \emptyset$, there exists $y_j\in L(v)\cap Y$, where $1\le j\le 2k+1$.
By Claim~\ref{cl11}, $N_H(y_j)\cap L(V(A))$ is a clique of $H$,
a contradiction to the fact that 
$\{x_j,x_{j+1}\}\subseteq N_H(y_j)$ and 
$\{x_j,x_{j+1}\}$ is an independent set of $H$.
Since $|L(v)\cap Y|\le k+1$, 
$\sigma(v,k,G,A,\cov)-|L(v)\cap Y|\ge 0$.

Claim~\ref{cl12} holds.
\claimend

We are now going to complete the proof.
Let 
$A_0=\{a\in A: L(a)\cap X\ne \emptyset\}$.
As $X$ is an independent set of $H$ with $|X|=2k+2$, 
we have $|A_0|=2k+2$.
Let $A'=A\setminus A_0$ and $\cov'=(L', H')$ be a cover of $G\vee A'$, where
$L'(u)=L(u)\setminus N_H(X)$ 
for all $u\in V(G\vee A')$
and $H'$ is the subgraph of $H$ induced by 
$\bigcup_{u\in V(G\vee A')}L'(u)$.

Consider the tuple $(k, G, A', \cov')$.
For each $u\in V(G\vee A')$, $|L'(u)|\ge |L(u)|-(\varn)$.
Thus, for each $a\in V(A')$, 
\eqn{th1-be17}
{
	|L'(a)|\ge |L(a)|-(\varn)
	\ge |A|+\chi(G)-(\varn)
	=|A'|+\chi(G).
}
By Claim~\ref{cl12}, for each $v\in V(G)$, $|L'(v)|\ge |L(v)|-(\varn-|L(v)\cap Y|)$.
Consequently, we have
\eqn{th1-e17}
{
	\sigma(v, k, G, A', \cov')&=&\max\{d_G(v)-|L'(v)|+|A'|+k,0\}
	\nonumber\\
	&\le& \max\{d_G(v)-(|L(v)|-(\varn-|L(v)\cap Y|))+|A|-(\varn)+k,0\}
	\nonumber\\
	&=&\sigma(v, k, G, A, \cov)-|L(v)\cap Y|,
}
where the last equality is due to the fact $\sigma(v, k, G, A, \cov) -|L(v)\cap Y|\ge 0$
by Claim~\ref{cl12}. 
Thus 
\eqn{th1-e18}
{
	\sigma(k, G, A', \cov')
	&\le& \sigma(k, G, A, \cov) - \sum_{v\in V(G)}|L(v)\cap Y|
	\nonumber \\
&=&\sigma(k, G, A, \cov) -(\vard)
}
and 
\eqn{th1-e19}
{
	|A'|&=&|A|-(\varn)
	\nonumber\\
	&\ge &\var\sigma(k, G, A, \cov)-(\varn)
	\nonumber\\
	&\ge &\var(\sigma(k, G, A', \cov')+\vard)-(\varn)
	\nonumber\\
	&=& \var\sigma(k, G, A', \cov').
}
By the assumption on the minimality of $|A|$, $G\vee A'$ is $\cov'$-colorable, 
whereas the union of $X$ and 
any $\cov'$-coloring of $G\vee A'$ is an $\cov$-coloring of $G\vee A$, a contradiction.
\claimend

Hence Theorem~\ref{th1} is proven.
\proofend

\section{
Proofs of Theorems~\ref{th1-5},~\ref{th1-3} and~\ref{th1-4}
\label{sec3}
}

We first prove the following lemma.

\begin{lemm}\label{le3-2}
$\chi_{DP}(C_4\vee K_q)=\chi(C_4\vee K_q)$ for all $q\ge 1$.
\end{lemm}

\proof
It suffices to show that $\chi_{DP}(C_4\vee K_1)=\chi(C_4\vee K_1)=3$.

Assume $\chi_{DP}(C_4\vee K_1)\ge 4$. 
Then there exists a 3-fold cover $\cov=(L,H)$ of $C_4\vee K_1$, 
such that $C_4\vee K_1$ has no $\cov$-coloring. 

For convenience purposes, let $G$ denote $C_4\vee K_1$.
Let $C_4$ be the cycle $x-y-u-v-x$ and $w$ be the vertex in $G$ 
adjacent to all vertices in $C_4$. 

For each vertex $z\in V(G)$, let 
$L(z)=\{z_1,z_2,z_3\}$. Thus, $L(x)=\{x_1,x_2,x_3\}$.

\setcounter{clm}{0}

\Clm{cl11-1}
{ 
For each edge $st\in E(G)$, 
$E_H(L(s),L(t))$ is a perfect matching.
}

\proof
Suppose that $st$ is an edge in $G$ 
such that $E_H(L(s),L(t))$ is not a perfect matching. 
Note that  $\{s,t\}\cap V(C_4)\ne \emptyset$. 
Without loss of generality, let $s$ be vertex $x$.

As $E_H(L(x),L(t))$ is not a perfect matching, 
$N_H(t_1)\cap L(x)=\emptyset$ for some $t_1\in L(t)$. 
Since $G\setminus \{x\}$ is a chordal graph, 
$G\setminus \{x\}$ has a perfect elimination ordering starting from $t$.
Thus, as $\chi(G\setminus \{x\})=3$ and $\cov$ is $3$-fold, 
$H\setminus L(x)$ has an independent set $I$ with $|I|=4$ and $t_1\in I$.
Since $d_G(x)=3$ and $N_H(t_1)\cap L(x)=\emptyset$, 
$I\cup \{x_i\}$ is an independent set of $H$ for some $x_i\in L(x)$,
implying that $G$ is $\cov$-colorable,
a contradiction.

Thus Claim~\ref{cl11-1} holds.
\claimend

Let $H^*$ be the spanning subgraph of $H$ obtained from $H$ 
by deleting all edges in the set $\{ab:a,b\in L(z),z\in V(G)\}$.

\Clm{cl11-2}
{ For any $z\in V(G)\setminus \{w\}$, 
$H^*\setminus L(z)$ has no $4$-cycles.
}

\proof
Assume Claim 2 fails. 

Suppose that $x_i-y_j-w_s-v_t-x_i$ is a cycle in $H^*\setminus L(u)$. 
Since $E_H(L(y),L(v))=\emptyset$,
$H\setminus (L(x)\cup L(w))$ has an independent set $I$ 
of size $3$ with $y_j, v_t\in I$.

As $x_i,w_s\in N_H(y_j)\cap N_H(v_t)$, 
$I$ can be extended to an independent set $I'$ of $H$ 
with $|I'|=5$ by choosing a suitable vertex in $L(w)$
and then a  suitable vertex in $L(x)$,
implying that $G$ is $\cov$-colorable, 
a contradiction.

Hence Claim~\ref{cl11-2} holds.
\claimend

\Clm{cl11-3}
{ 
For each vertex $b\in V(H)\setminus L(w)$, 
$H^*[N_{H^*}(b)]$ is a path.  
}

\proof Suppose the claim fails. 
Without loss of generality, suppose that 
$H^*[N_{H^*}(x_1)]$ is not a path.
By Claim~\ref{cl11-1}, let $N_{H^*}(x_1)=\{y_j,w_s,v_t\}$.
Then, $\{y_j,v_t\}\not\subseteq N_H(w_s)$.
Assume that $y_j\notin N_H(w_s)$.

Then, $\{y_j,w_s\}$ can be extended to an independent set $I$ of $H$ 
with $|I|=5$ 
by choosing a suitable vertex in $L(u)$,
then a suitable vertex in $L(v)$,
and finally a suitable vertex in $L(x)$.
Thus,  $G$ is $\cov$-colorable, 
a contradiction.
\claimend

But, Claims~\ref{cl11-1} and~\ref{cl11-3} together imply that $H^*\setminus L(u)$ contains $4$-cycles,
a contradiction to Claim~\ref{cl11-2}.

Hence $\chi_{DP}(C_4\vee K_1)=3$ and the result holds.
\proofend

We will also apply the following result 
in the proof of Theorem~\ref{th1-5}.

\begin{lemm}[\cite{Dong2004}] \label{le3-1}
Let $G$ be a graph with $n$ vertices, where $n\ge 3$.
If $\chi(G)=n-2$, then either $G\cong C_5\vee K_{n-5}$ or 
$G\setminus \{u,v\}\cong K_{n-2}$ for some vertices $u$ and $v$.
\end{lemm}

Now we are going to complete the proofs of Theorems~\ref{th1-5},
~\ref{th1-3} and~\ref{th1-4}.

\noindent {\it Proof of Theorem~\ref{th1-5}}: 
(i) To show that $Z_{DP}(C_4)=1$, 
by Lemma~\ref{le3-2},
it suffices to show that $\chi_{DP}(C_4)\ne 2$.

\begin{figure}[!ht]
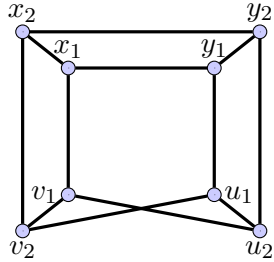

	\tikzp{1.2}
	{
		\foreach \place/\y in {	{(-0.8,0.7)/1},
			{(0.8,0.7)/2},{(0.8,-0.7)/3},{(-0.8,-0.7)/4},{(-1.3,1.1)/5}, {(1.3,1.1)/6},{(1.3,-1.1)/7},{(-1.3,-1.1)/8}}   
		\node[cblue] (b\y) at \place {};
		
		\filldraw[black] (b1) circle (0pt)node[anchor=south] {$x_1$};
		\filldraw[black] (b2) circle (0pt)node[anchor=south] {$y_1$};
		\filldraw[black] (b3) circle (0pt)node[anchor=west] {$u_1$};
		\filldraw[black] (b4) circle (0pt)node[anchor=east] {$v_1$};
		\filldraw[black] (b5) circle (0pt)node[anchor=south] {$x_2$};
		\filldraw[black] (b6) circle (0pt)node[anchor=south] {$y_2$};
		\filldraw[black] (b7) circle (0pt)node[anchor=north] {$u_2$};
		\filldraw[black] (b8) circle (0pt)node[anchor=north] {$v_2$};

		\draw[black, very thick] (b1) -- (b2) -- (b3)--(b8) -- (b5) -- (b6) -- (b7)--(b4)--(b1);  
		\draw[black, very thick] (b1) -- (b5);  
		\draw[black, very thick] (b2) -- (b6);  
		\draw[black, very thick] (b3) -- (b7);  
		\draw[black, very thick] (b4) -- (b8);  
	}
	\caption{Graph $H_0$}
	\label{fig3-2-1}
\end{figure}

Consider the $2$-fold cover $\cov=(L,H_0)$, where 
$H_0$ is as shown in Figure~\ref{fig3-2-1}.
Assume that $I$ is an independent set of $H$ with $|I|=4$.
If $x_i\in I$, where $i\in \{1,2\}$, 
then $y_{3-i},v_{3-i}\in I$, implying that 
$\{u_1,u_2\}\cap I=\emptyset$, a contradiction.
Thus, $\chi_{DP}(C_4)\ne 2$.

(ii) 
Assume that $\chi_{DP}(G)>\chi(G)=n-2$.
By Lemma~\ref{le3-1}, 
 either $G\cong C_5\vee K_{n-5}$ or 
$G\setminus \{u,v\}\cong K_{n-2}$ for some vertices $u$ and $v$.

As $\chi_{DP}(G)>\chi(G)$, by (\ref{1-1}), we have 
$col(G)>\chi(G)=n-2$.
If  $G\cong C_5\vee K_{n-5}$, it can be verified that 
$col(G)=\chi(G)=n-2$, 
a contradiction.

Thus, $G\setminus \{u,v\}\cong K_{n-2}$ for some vertices $u$ and $v$.
If $uv\notin E(G)$, then both $u$ and $v$ are simplicial,
implying that $G$ is chordal and $\chi_{DP}(G)=\chi(G)$ by Lemma~\ref{le1-1},
a contradiction again.
Hence $uv\in E(G)$.

If $V(G)\setminus \{u,v\}$ is a subset of $N_G(u)$ or $N_G(v)$,
then $\chi(G)\ge n-1$, a contradiction. 
Thus, $d(u)\le n-2$ and $d(v)\le n-2$.

If $d(u)\le n-3$, as $uv\in E(G)$ and $d(v)\le n-2$, 
it can be verified that $col(G)=n-2=\chi(G)$,
a contradiction. 
Thus, $d(u)=n-2$. Similarly, $d(v)=n-2$.
It implies that $N_G(u)=V(G)\setminus \{u,u'\}$ 
and $N_G(v)=V(G)\setminus \{v,v'\}$
for some vertices $u'$ and $v'$ in $V(G)\setminus \{u,v\}$.
If  $u'=v'$, then $G\setminus u'\cong K_{n-1}$, a contradiction. 
Thus, $u'\ne v'$.

Clearly, $n\ge 4$ and the above conclusion implies that 
$G$ is isomorphic to the graph obtained from $K_{n}$ by removing 
two edges which have no common end,
in other words, $G\cong C_4\vee K_{n-4}$.

By Lemma~\ref{le3-2}, $\chi_{DP}(C_4\vee K_s)=\chi(C_4\vee K_s)$
for each $s\ge 1$.
Hence Theorem~\ref{th1-5} holds.
\proofend

By Lemma~\ref{le1-1}, Theorems~\ref{th1} and~\ref{th1-5}, we give the proof of Theorem~\ref{th1-3} below.

\noindent {\it Proof of Theorem~\ref{th1-3}}: 

If $k=1$ or $k\ge n-1$, then 
$G$ is chordal,
and thus $\chi_{DP}(G)=\chi(G)$ by Lemma~\ref{le1-1},
implying that $Z_{DP}(G)=0$.

If $G=C_4$, then $Z_{DP}(G)=1$ by Theorem~\ref{th1-5}.
If $k=n-2$ but $G\not\cong C_4$,
then $Z_{DP}(G)=0$ by Theorem~\ref{th1-5} again.

Now assume that $2\le k\le n-3$. 
Let $r=\lceil \frac{4(k+1)m}{2k+1} \rceil$ 
and $A$ be a complete graph with $r$ vertices.
For any ($k+r$)-fold cover $\cov$ of $G\vee A$, 
we have $|L(u)|=|A|+\chi(G)$ for all $u\in V(G\vee A)$ and 
\equ{co-e1}
{
	|A|=r
	\ge \var \sum_{v\in V(G)}d_G(v) 
	= \var\sum_{v\in V(G)}\max\{d_G(v)-|L(v)|+|A|+k,0\}.
}
By Theorem~\ref{th1}, 
$G\vee A$ is $\cov$-colorable and 
$\chi_{DP}(G\vee A)\le k+r=\chi(G\vee A)$,
implying that $\chi_{DP}(G\vee A)=\chi(G\vee A)$,
i.e., $\chi_{DP}(G\vee K_r)=\chi(G\vee K_r)$.
Hence $Z_{DP}(G)\le r=
\lceil \frac{4(k+1)m}{2k+1} \rceil 
\le \lceil 2.4m\rceil$.
\proofend

Finally, we prove Theorem~\ref{th1-4} below.

\noindent {\it Proof of Theorem~\ref{th1-4}}:
Let $G$ be any graph of order $n$
and $k=\chi(G)$.

By Lemma~\ref{le1-1} and Theorem~\ref{th1-5}, if $k=1$ or $k\ge n-2$ and $G\not\cong C_4$, then $Z_{DP}(G)=0< n^2-(n+3)/2$ when $n\ge 2$;
if $G=C_4$, then $Z_{DP}(G)=1<4^2-(4+3)/2$.

Now assume that $2\le k\le n-3$. 
It can be verified easily that $|E(G)|\le (n/k)^2{k\choose 2}$. 
Thus, by Theorem~\ref{th1-3}, 
\eqn{co-e2}
{
	Z_{DP}(G) 
	&\le & 1+\frac{4(k+1)}{2k+1}\times (n/k)^2{k\choose 2} 
	\nonumber \\
	&=&1+\frac{2(k^2-1)}{k(2k+1)}n^2
	\nonumber \\
	&=&n^2+1-\frac{(k+2)n^2}{2k^2+k}
	\nonumber \\
	&\le &
	n^2+1-\frac{(n-1)n^2}{2(n-3)^2+(n-3)}
	\nonumber \\
	&= &
	n^2+1-\frac{n}{2}-\frac{3n(3n-5)}{2(n-3)(2n-5)}
	\nonumber \\
	&< &
	n^2-\frac{n}{2}-\frac{5}{4}. 
}
As $Z_{DP}(G)$ is an integer, we have $Z_{DP}(G)\le n^2-\frac{n+3}2$.
Hence $Z_{DP}(n)\le n^2-\frac{n+3}2$.
\proofend

\end{document}